\documentclass[a4paper,leqno]{amsart}

\usepackage{amsmath,amssymb,amsthm}

\def\C{{\mathbb C}}
\def\R{{\mathbb R}}
\def\N{{\mathbb N}}
\def\T{{\mathbb T}}

\def\O{\mathcal O}

\def\virgp{\raise 2pt\hbox{,}}
\def\bv{{\bf v}}
\def\bu{{\bf u}}
\def\bw{{\bf w}}
\def\({\left(}
\def\){\right)}
\def\<{\left\langle}
\def\>{\right\rangle}
\def\le{\leqslant}
\def\ge{\geqslant}

\DeclareMathOperator{\RE}{Re}
\DeclareMathOperator{\IM}{Im}
\def\defn{\mathrel{:=}}

\DeclareMathOperator{\DIV}{\mathrm{div}}

\def\d{{\partial}}

\def\eps{\varepsilon}
\def\l{\lambda}

\theoremstyle{plain}
\newtheorem{theorem}{Theorem}[section]
\newtheorem{lemma}[theorem]{Lemma}
\newtheorem{corollary}[theorem]{Corollary}
\newtheorem{proposition}[theorem]{Proposition}
\newtheorem*{hyp}{Assumptions}

\theoremstyle{definition}

\theoremstyle{remark}
\newtheorem{remark}[theorem]{Remark}
\newtheorem*{remark*}{Remark}

\numberwithin{equation}{section}

\begin{document}

\title[WKB analysis for the Gross--Pitaevskii equation]{WKB analysis
  for the Gross--Pitaevskii equation with non-trivial boundary
  conditions at infinity}
 \author[T. Alazard]{Thomas Alazard}
\address{CNRS \& Universit\'e Paris-Sud\\ Math\'ematiques UMR CNRS
 8628 \\ B\^at. 425\\ 91405
  Orsay cedex\\ France}
\email{Thomas.Alazard@math.cnrs.fr}
\author[R. Carles]{R{\'e}mi Carles}
\address{CNRS \& Universit\'e Montpellier~2\\ Math\'ematiques
  UMR CNRS 5149\\ CC 051\\ Place Eug\`ene Bataillon\\ 34095
  Montpellier cedex 5\\ France}
\email{Remi.Carles@math.cnrs.fr}
\begin{abstract}
We consider the semi-classical limit for the Gross--Pitaevskii
equation. In order to consider
non-trivial boundary conditions at infinity, we work in Zhidkov spaces
rather than in Sobolev spaces. For the usual cubic nonlinearity, we
obtain a point-wise description of the wave function as the Planck
constant goes to zero, so long as no singularity appears in the limit
system. For a cubic-quintic nonlinearity, we show that working with
analytic data may be necessary and sufficient to obtain a similar result.
\end{abstract}
\maketitle

\section{Introduction}
\label{sec:intro}

We study the semi-classical limit $\hbar \to 0$ for the
Gross--Pitaevskii equation
\begin{equation*}
  i\hbar \d_t u + \frac{\hbar^2}{2m}\Delta u = V u + f\( |u|^2\)u,
\end{equation*}
where $x\in \R^n$. In the case of Bose--Einstein condensation (BEC), the
external potential $V=V(t,x)$ models an
external trap, and the
nonlinearity $f$ describes the nonlinear
interactions of the particles  (see
e.g. \cite{DGPS,PiSt,JosserandPomeau}). We consider two types of
nonlinearity $f$ (after renormalization):
\begin{itemize}
\item Cubic nonlinearity: $f(|u|^2)u=\(|u|^2-1\)u$.
\item Cubic-quintic nonlinearity: $f(|u|^2)u=\(|u|^4+\lambda
  |u|^2\)u$, $\lambda \in \R$.
\end{itemize}
The cubic nonlinearity is certainly the most commonly used model in
BEC. The
defocusing nonlinearity corresponds to a positive scattering length,
as in the case of $\,^{87}$Rb, $\,^{23}$Na and $\,^1$H. Note that this
model is also used in superfluid theory. See
e.g. \cite{DGPS,PiSt,JosserandPomeau} and references therein. The
cubic-quintic nonlinearity, which is mostly used as an envelope
equation in optics, is also considered in BEC for alkalimetal gases (see
e.g. \cite{JPhysB,PhysRevA63,PhysRevE}), in which case $\lambda
<0$. The cubic term corresponds to
a negative scattering length, and the quintic term to a repulsive
three-body elastic interaction. We also consider the case $\l>0$
(positive scattering length).

\subsection{Cubic nonlinearity}

Up to rescaling the Planck constant, we consider the limit $\eps\to 0$
for:
\begin{align}
  i\eps \d_t u^\eps +\frac{\eps^2}{2}\Delta u^\eps &= V u^\eps+\(
  |u^\eps|^2-1\)u^\eps,\quad x\in \R^n,\ n\ge 1,\label{eq:gp}\\
u^\eps(0,x)&=a_0^\eps(x)e^{i\phi_0(x)/\eps}.\label{eq:ci}
\end{align}
Our initial data  do not necessarily decay to zero at
infinity. Typically, we do not assume $a_0^\eps\in
L^2(\R^n)$ (see Theorem~\ref{theo:main} below). Recently, the Cauchy
problem  \cite{Gallo,PG05} and the
semi-classical limit \cite{LinZhang} for \eqref{eq:gp} with $V\equiv
0$ have been studied more systematically.
When the external potential $V$ is zero, $V\equiv 0$, the Hamiltonian
structure yields, at least formally:
\begin{equation*}
  \frac{d}{dt}\(\|\eps \nabla u^\eps(t)\|_{L^2}^2 +
  \left\||u^\eps(t)|^2-1\right\|_{L^2}^2  \)=0.
\end{equation*}
In this case, a natural
space to study the Cauchy problem associated to \eqref{eq:gp}  is the
energy space (see e.g. \cite{BethuelSaut,PG05} and references therein)
\begin{equation*}
  E=\{ u \in H^1_{\rm loc}(\R^n)\ ;\ \nabla u\in L^2(\R^n), \ |u|^2-1
  \in L^2(\R^n)\}.
\end{equation*}
For this quantity to be well defined, one cannot assume that $u^\eps$ is
in $L^2(\R^n)$; morally, the modulus of $u^\eps$
goes to one at infinity. To study solutions which are bounded, but not
in $L^2(\R^n)$, P.~E.~Zhidkov introduced in the
one-dimensional case in \cite{Zhidkov} (see also \cite{ZhidkovLNM}):
\begin{equation}\label{eq:zhidkov}
  X^s(\R^n) = \{  u \in L^\infty(\R^n)\ ;\ \nabla u \in
  H^{s-1}(\R^n)\},\quad s>n/2.
\end{equation}
We also denote $X^\infty := \cap_{s>n/2}X^s$.
The study of these spaces was
generalized in the multidimensional case by C.~Gallo
\cite{Gallo}. They make it possible to consider solutions to
\eqref{eq:gp} whose modulus has a non-zero limit as $|x|\to
\infty$, but not necessarily satisfying $|u^\eps(t,\cdot)|^2 -1\in
L^2(\R^n)$. We shall also use these spaces.
\smallbreak

Recently, P.~G\'erard \cite{PG05} has solved the Cauchy problem for the
Gross--Pitaevskii equation in the more natural space $E$, in space
dimensions two and three. The main novelty consists in working
with distances instead of norms, in order to apply a fixed point
argument in $E$. In particular, the constraint  $|u^\eps(t,\cdot)|^2 -1\in
L^2(\R^n)$ is satisfied.
\smallbreak

To our knowledge, if the initial data do not vanish at infinity, the
introduction of an (unbounded) external
potential in Gross--Pitaevskii equation has no physical
motivation. Note also that if $V$ is an harmonic potential, then the
formal Hamiltonian corresponding to \eqref{eq:gp} is necessarily
infinite (see \S~\ref{sec:hamilton}).
On the other hand, introducing
a quadratic external potential or considering a quadratic initial phase
$\phi_0$ makes no difference in our analysis.
The model
\eqref{eq:gp}--\eqref{eq:ci} with $V\equiv 0$ and $\phi_0$ quadratic is
certainly more physically relevant, and does not seem to enter into
the framework of the previous mathematical studies. 
Another motivation to 
introduce this external potential stems from the 
study of the semi-classical limit of the Schr\"odinger--Poisson system, where
$|u^\eps|^2-1$ is replaced with 
$V_{\rm p}^\eps$ given by 
$\Delta V_{\rm p}^\eps = q\(|u^\eps|^2-c\)$. 
This models appears in the semi-conductor theory 
where the real number $q$ models a charge, which we may
take equal to one here, and the function $c=c(x)$ models a doping
profile, which we may take to be $c\equiv 1$. 
As in \cite{AC-SP}, we will prove that if $V$ grows
quadratically in space, then if $|u^\eps(t=0,\cdot)|^2-1\in
L^2(\R^n)$, one must not expect $|u^\eps(t,\cdot)|^2-1\in
L^2(\R^n)$ for $t>0$.

\begin{hyp}\label{hyp:geom}
  We assume that the potential and the initial phase are of the form:
  \begin{itemize}
  \item $V\in C^\infty(\R_t\times \R^n_x)$, and $V=V_{\rm quad}+V_{\rm
  lin}$, where $V_{\rm quad}(t,x)=\,^{t}xM(t)x$ is a quadratic form,
  with $M(t)\in {\mathcal S}_{n}(\R)$ a symmetric $n\times n$ matrix,
  depending smoothly on $t$, and $\nabla V_{\rm
  lin} \in C^\infty(\R_t;X^s)$ for all $s>n/2$.
  \item $\phi_0\in C^\infty( \R^n)$, and  $\phi_0=\phi_{\rm quad}+\phi_{\rm
  lin}$, where $\phi_{\rm quad}(x)=\,^{t}xQ_0x$ is a quadratic form,
  with $Q_0$ a symmetric matrix in ${\mathcal M}_{n\times n}(\R)$, and
  $\nabla \phi_{\rm lin} \in X^\infty$.
  \end{itemize}
\end{hyp}
Note that our assumptions include the case where $V_{\rm lin}$ and
$\phi_{\rm lin}$ are linear in $x$. In general, these functions are
sub-linear in $x$, since their gradient is bounded.
\begin{lemma}\label{lem:hj}
  There exist $T>0$ and a unique
  solution $\phi_{\rm eik}\in C^\infty([0,T]\times\R^n)$ to:
  \begin{equation}
    \label{eq:eik}
    \partial_t \phi_{\rm eik} +\frac{1}{2}|\nabla_x \phi_{\rm eik}|^2
    +V_{\rm quad}=0\quad ;\quad
    \phi_{{\rm eik} \mid t=0}=\phi_{\rm quad}\, .
  \end{equation}
Moreover, $\phi_{\rm eik}$ is a quadratic form in
$x$:
\begin{equation}\label{eq:formphieik}
  \phi_{\rm eik}(t,x) = \,^{t}xQ(t)x ,
\end{equation}
where $Q(t)\in {\mathcal S}_{n}(\R)$ is a smooth function of $t$.
\end{lemma}
\begin{proof}
  Existence and uniqueness  follow from \cite[Lemma~1]{CaBKW}. To
  prove that $\phi_{\rm eik}$ is quadratic in $x$, seek $\phi_{\rm
  eik}$ of the form \eqref{eq:formphieik}. Then \eqref{eq:eik} is
  equivalent to the system of ordinary differential equations
  \begin{align*}
    \dot Q(t) +2Q(t)^2 + M(t) =0\quad ;\quad Q(0)=Q_0.
  \end{align*}
The lemma then follows from Cauchy--Lipschitz Theorem.
\end{proof}
\begin{remark}\label{rema:transport}
As in \cite{AC-SP}, we shall use the following geometrical interpretation of the above
  lemma. The time $T$ is such that for $t\in [0,T]$, the map given by
\begin{equation*}
  \partial_t x(t,y) = \nabla_x \phi_{\rm eik}
  \left(t,x(t,y)\right)= Q(t)x(t,y)\quad ; \quad x(0,y)=y,
\end{equation*}
defines a global diffeomorphism on $\R^n$.
Therefore, the characteristics associated to
  the operator $\d_t +\nabla \phi_{\rm eik}\cdot \nabla$ do not meet
  for $t\in [0,T]$, and this operator is a smooth transport
  operator:
  \begin{equation*}
   \(\d_t +\nabla \phi_{\rm eik}\cdot \nabla\)(f(t,y))= \d_t
   f\(t,x(t,y)\).
  \end{equation*}
Note that if $Q(t)$ and its anti-derivative commute, then we
  have
  \begin{equation*}
    x(t,y)=\exp\( \int_0^t Q(\tau)d\tau \)y.
  \end{equation*}
\end{remark}

\begin{theorem}\label{theo:main}
 Suppose that there exist $a_0, a_1\in X^\infty$ such that:
  \begin{equation}\label{eq:DACI}
    \left\|a_0^\eps -a_0 -\eps a_1\right\|_{X^s}
    =o(\eps),\quad \forall s>n/2.
  \end{equation}
There exist $T_*\in ]0,T]$ independent
  of $\eps \in ]0,1]$, and a unique solution $u^\eps \in
  C^\infty\cap L^\infty([0,T_*]\times \R^n)$  to
  \eqref{eq:gp}--\eqref{eq:ci}. Moreover, there exist $a,\phi\in
  C^\infty([0,T_*]\times\R^n)$ with
  $a,\nabla \phi \in C([0,T_*];X^s)$ for all $s>n/2$, such that:
  \begin{equation}\label{eq:BKWVtpetit}
    \limsup_{\eps \to 0}\left\|
    u^\eps (t,\cdot) - a(t,\cdot)
    e^{i(\phi(t,\cdot)+\phi_{\rm eik}(t,\cdot))/\eps}\right\|_{
    L^\infty(\R^n)}=\O(t)\quad \text{as }t\to 0.
  \end{equation}
The functions $a$ and $\phi$ depend nonlinearly on
$\phi_0$ and
$a_0$ (see \eqref{eq:systexact0} below). There exists
$\phi^{(1)}\in L^\infty([0,T_*]\times
    \R^n)$, real-valued, with $\nabla \phi^{(1)}\in
  C([0,T_*];X^s)$ for all $s>n/2$, such that:
  \begin{equation}\label{eq:BKWVtgrand}
    \limsup_{\eps \to 0}\left\|
    u^\eps - ae^{i\phi^{(1)}}
    e^{i(\phi+\phi_{\rm eik})/\eps}\right\|_{L^\infty([0,T_*]\times
    \R^n)}=0.
  \end{equation}
The modulation $\phi^{(1)}$ is a nonlinear function of $\phi_0$, $a_0$
and $a_1$ (see \eqref{eq:systlinear} below).
\end{theorem}
\begin{remark}
  Several applications of this general results are given, in
  \S\ref{sec:semi}, \S\ref{sec:modone} and \S\ref{sec:hamil}.
\end{remark}
\begin{remark}\label{rema:mieux}
If we assume moreover
\begin{equation*}
    \left\|a_0^\eps -a_0 -\eps a_1\right\|_{X^s}
    =\O(\eps^2),\quad \forall s>n/2,
  \end{equation*}
then the above error estimate can be improved:
\begin{equation*}
    \left\|
    u^\eps - ae^{i\phi^{(1)}}
    e^{i(\phi+\phi_{\rm eik})/\eps}\right\|_{L^\infty([0,T_*]\times
    \R^n)}=\O(\eps).
  \end{equation*}
\end{remark}
\begin{remark}
  The above result  and System~\eqref{eq:systlinear} below show that
  in general, it is necessary to know the initial amplitude $a_0^\eps$
  up to the order $o(\eps)$ to describe the leading order behavior of the
  \emph{wave function} $u^\eps$. It is not necessary to know
  $a_0^\eps$ with such precision to study the convergence of quadratic
  observables. See \S\ref{sec:hydro}. In particular, in 
  Theorem~\ref{theo:hydro}, we extend the result of \cite{LinZhang} to
  the three-dimensional case (on a bounded domain, or outside a
  bounded domain). 
\end{remark}
\begin{remark}
Most of the results that we present here remain valid in a
space-periodic setting, that is if we assume $x\in \T^n$.  In that
case, compactness arguments show that the proof of
Theorem~\ref{theo:main} remains valid when $V\in C^\infty(\R_t\times
\T^n_x)$ and $\phi_0\in C^\infty(\T^n)$. On the other hand, the
discussions in \S\ref{sec:hamilton} and \S\ref{sec:hamil}
become irrelevant on the torus. Finally, note that it is equivalent to
work in Sobolev spaces, since $X^s(\T^n)=H^s(\T^n)$ for
$s>n/2$.
\end{remark}

The analysis detailed in \S\ref{sec:general} and \S\ref{sec:semi}
shows that the formal part of \cite{CaARMA} can be justified in the present
framework. We shall only state a typical consequence of this approach:
\begin{corollary}[Instability]
  Let $n\ge 1$, $a_0 ,a_1\in C^\infty\cap
  X^\infty(\R^n)$, with $\RE(\overline a_0 a_1)\not \equiv 0$, and
  $\phi_0\in C^\infty(\R^n)$, with $\nabla \phi_0\in X^\infty$. Let
  $u^\eps$ and
$v^\eps$ solve the initial value problems:
\begin{align*}
i\eps \d_t u^\eps + \frac{\eps^2}{2}\Delta u^\eps &=
\left(|u^\eps|^2-1\right)u^\eps
\ ; \ u^\eps\big|_{t=0}= a_0e^{i\phi_0/\eps}\, .\\
i\eps \d_t v^\eps + \frac{\eps^2}{2}\Delta v^\eps &=
\left(|v^\eps|^2-1\right)v^\eps
\ ; \ v^\eps\big|_{t=0}= \(a_0+ \delta^\eps a_1\)e^{i\phi_0/\eps}\, ,
\end{align*}
where $\delta^\eps \to 0$.
Assume that there exists $N\in {\mathbb N}$ such that
$\delta^\eps/\eps^{1-\frac{1}{N}}\to +\infty$.
Then we can find $t^\eps \to 0$ such that
$\displaystyle \liminf_{\eps \to 0}\left\| u^\eps(t^\eps) - v^\eps(t^\eps)
\right\|_{L^\infty}>0$. In particular,
\begin{equation*}
\liminf_{\eps \to 0}\frac{\left\| u^\eps - v^\eps
  \right\|_{L^\infty([0,t^\eps]\times \R^n)}}{\left\| u^\eps_{\mid
  t=0} - v^\eps_{\mid t=0}
  \right\|_{L^\infty(\R^n)}}= +\infty.
\end{equation*}
\end{corollary}
\begin{remark}
  Note that if $\phi_0\equiv 0$, then we also have:
\begin{equation*}
\liminf_{\eps \to 0}\frac{\left\| u^\eps - v^\eps
  \right\|_{L^\infty([0,t^\eps]\times \R^n)}}{\left\| u^\eps_{\mid
  t=0} - v^\eps_{\mid t=0}
  \right\|_{X^s}}= +\infty,\quad \forall s>n/2.
\end{equation*}
This shows that the instability mechanism is not due to regularity
issues. It is due to the fact that \eqref{eq:gp} is super-critical as
far as WKB analysis is concerned: the small initial perturbation (of
order $\delta^\eps$)
yields a high-frequency perturbation of the evolution (a
multiplicative factor of the form
$e^{-2it \delta^\eps\RE(\overline a_0 a_1) /\eps}$).
\end{remark}
\subsection{Cubic-quintic nonlinearity}
Denote $f_\l (y) = y^2 +\l y$. We now consider
\begin{equation}\label{eq:cubicquintic}
\left\{
  \begin{aligned}
    i\eps \d_t u^\eps +\frac{\eps^2}{2}\Delta u^\eps &= f_\l \(
  |u^\eps|^2\)u^\eps,\quad x\in \R^n,\ n\ge 1,\\
u^\eps(0,x)&=a_0^\eps(x)e^{i\phi_0(x)/\eps}.
  \end{aligned}
  \right.
\end{equation}
Note that in \eqref{eq:cubicquintic}, we assume that there is no
external potential, $V=0$. We also assume that there is no initial
quadratic oscillation: $\phi_0\in C^\infty(\R^n;\R)$, with $\nabla
\phi_0\in X^\infty$. The case $\l>0$, $V\not =0$, with $a_0^\eps\in
H^\infty$, is contained in \cite{CaBKW}. We assume $V_{\rm quad}=0$
here in order to consider non-zero boundary conditions at infinity. We
also assume $V_{\rm lin}=0$ for simplicity only.
\smallbreak

Plugging an approximate solution of the form $u^\eps\thickapprox a
e^{i\phi/\eps}$, with $a$ and $\phi$ independent of $\eps$, and
passing to the limit $\eps\to 0$ as in
\cite{GasserLinMarkowich,LinZhang}, we find formally that
$(\rho,v)\defn (|a|^2,\nabla \phi)$ solves:
\begin{equation}\label{eq:eulerhypell}
  \left\{
    \begin{aligned}
      &\d_t \rho +\DIV \(\rho v\) = 0.\\
      & \d_t v + v\cdot \nabla v + \nabla \(f_\l(\rho)\)=0.
    \end{aligned}
\right.
\end{equation}
If $\l>0$, then the problem is hyperbolic. Essentially, the result of
Theorem~\ref{theo:main} remains valid. When $\l<0$,
the above problem is hyperbolic for $\rho>|\l|/2$ and elliptic for
$\rho<|\l|/2$. This feature is reminiscent of Euler equations of gas
dynamics in Lagrangian coordinates:
\begin{equation}\label{eq:eulerp}
  \left\{
    \begin{aligned}
      &\d_t u +\d_x v= 0.\\
      & \d_t v + \d_x \(p(u)\)=0.
    \end{aligned}
\right.
\end{equation}
As recalled in \cite{GuyCauchy}, a typical mathematical example for
van der Waals state laws is given by $p(u)=(u^2-1)u$. The problem is
hyperbolic 
if $u>1/\sqrt 3$, and elliptic if $u<1/\sqrt 3$. Hadamard's argument
implies that the only
reasonable framework to study \eqref{eq:eulerhypell} or
\eqref{eq:eulerp} is that of analytic functions (see
\cite{GuyCauchy}). In this case,
we refer to the approach of \cite{PGX93,ThomannAnalytic}. More details
are given in \S\ref{sec:quintic}. When the elliptic region for
\eqref{eq:eulerhypell} is avoided, then essentially,
Theorem~\ref{theo:main} remains valid:
\begin{theorem}\label{theo:cubicquintic}
 Suppose that there exist $a_0, a_1\in X^\infty$ such that:
  \begin{equation*}
    \left\|a_0^\eps -a_0 -\eps a_1\right\|_{X^s}
    =o(\eps),\quad \forall s>n/2.
  \end{equation*}
Assume moreover that $\phi_0\in C^\infty(\R^n;\R)$ with $\nabla
\phi_0\in X^\infty$, and:
\begin{itemize}
\item Either $\l>0$,
\item Or $\l<0$ and there exists $\delta>0$ such that $|a_0(x)|^2\ge
  \delta+ \frac{|\l|}{2}$, $\forall x\in \R^n$.
\end{itemize}
Then there exist $\eps_*,T_*>0$, and a unique solution $u^\eps \in
  C^\infty\cap L^\infty([0,T_*]\times \R^n)$  to
  \eqref{eq:cubicquintic} for all $\eps \in
  ]0,\eps_*]$. Moreover, there exist $a,\phi\in
  C^\infty([0,T_*]\times\R^n)$ with
  $a,\nabla \phi \in C([0,T_*];X^s)$ for all $s>n/2$, such that:
  \begin{equation*}
    \limsup_{\eps \to 0}\left\|
    u^\eps (t,\cdot) - a(t,\cdot)
    e^{i\phi(t,\cdot)/\eps}\right\|_{
    L^\infty(\R^n)}=\O(t)\quad \text{as }t\to 0.
  \end{equation*}
There exists
$\phi^{(1)}\in L^\infty([0,T_*]\times
    \R^n)$, real-valued, with $\nabla \phi^{(1)}\in
  C([0,T_*];X^s)$ for all $s>n/2$, such that:
  \begin{equation*}
    \limsup_{\eps \to 0}\left\|
    u^\eps - ae^{i\phi^{(1)}}
    e^{i\phi/\eps}\right\|_{L^\infty([0,T_*]\times
    \R^n)}=0.
  \end{equation*}
\end{theorem}

\subsection{Structure of the paper}

In \S\ref{sec:general}, we construct the solution $u^\eps$ as $u^\eps
= a^\eps e^{i\Phi^\eps/\eps}$, where $a^\eps$ is complex-valued and
$\Phi^\eps$ is real-valued. This yields the existence part of
Theorems~\ref{theo:main} and \ref{theo:cubicquintic}. The proof of
these theorems is completed in \S\ref{sec:semi}, where the limit of
$(a^\eps,\Phi^\eps)$ as $\eps$ goes to zero is studied. We give three
examples of applications of Theorem~\ref{theo:main} in
\S\ref{sec:modone}, in the case $\phi_{\rm eik}=0$. In
\S\ref{sec:hamil}, we study the time evolution
of a non-trivial boundary condition at infinity when $\phi_{\rm
  eik}\not =0$. In \S\ref{sec:hydro}, we investigate the limit of the
position and current densities. Finally, we explain why working in an
analytic setting is often necessary (and always sufficient) in the
case of the cubic-quintic nonlinearity.

\section{Construction of the solution}
\label{sec:general}

\subsection{Phase-amplitude representation: the case $\phi_{\rm eik}=V=0$}

When $V$ and  $\phi_0$ are
identically zero, the existence and uniqueness part of
Theorem~\ref{theo:main} was established by C.~Gallo \cite{Gallo}. Note
however that with our scaling, the fact that $T_*$ is independent of
$\eps\in ]0,1]$ does not follow from \cite{Gallo}. Since the
approach in Zhidkov spaces is rather similar to the one in Sobolev
spaces, we shall essentially explain the new aspects of
the proof. To treat both cubic and cubic-quintic nonlinearities,
consider the general equation
\begin{equation}\label{eq:general}
\left\{
  \begin{aligned}
    i\eps \d_t u^\eps +\frac{\eps^2}{2}\Delta u^\eps &= f \(
  |u^\eps|^2\)u^\eps,\quad x\in \R^n,\ n\ge 1,\\
u^\eps(0,x)&=a_0^\eps(x)e^{i\phi_0(x)/\eps},
  \end{aligned}
  \right.
\end{equation}
where $f\in C^\infty (\R_+;\R)$. We keep
the hierarchy introduced by E.~Grenier \cite{Grenier98}: seek
$u^\eps=a^\eps e^{i\Phi^\eps}$, where $a^\eps$ is complex-valued, and
$\Phi^\eps$ is real-valued. We impose
\begin{equation}\label{eq:systemmanuel}
  \left\{
\begin{aligned}
    \partial_t \Phi^\eps +\frac{1}{2}\left|\nabla
    \Phi^\eps\right|^2 +
    f\(|a^\eps|^2\)= 0\quad &; \quad
    \Phi^\eps\big|_{t=0}=\phi_0,\\
\partial_t a^\eps +\nabla \Phi^\eps \cdot \nabla
    a^\eps +\frac{1}{2}a^\eps
\Delta \Phi^\eps  = i\frac{\eps}{2}\Delta
    a^\eps\quad & ;\quad
a^\eps\big|_{t=0}= a^\eps_0\, .
\end{aligned}
\right.
\end{equation}
As an intermediary unknown function, introduce the ``velocity''
$v^\eps = \nabla \Phi^\eps$. Separate real and
imaginary parts of $a^\eps$, $a^\eps =
a_1^\eps + ia_2^\eps$, and introduce:
\begin{equation*}
  \bu^\eps = \left(
    \begin{array}[l]{c}
       a_1^\eps \\
       a_2^\eps \\
       v^\eps_1 \\
      \vdots \\
       v^\eps_n
    \end{array}
\right)\ , \quad
\bu^\eps_0 = \left(
    \begin{array}[l]{c}
       \RE(a_0^\eps) \\
       \IM(a_0^\eps) \\
       \d_1 \phi_0 \\
      \vdots \\
       \d_n \phi_0
    \end{array}
\right)
\ , \quad
L = \left(
  \begin{array}[l]{ccccc}
   0  &-\Delta &0& \dots & 0   \\
   \Delta  & 0 &0& \dots & 0  \\
   0& 0 &&0_{n\times n}& \\
   \end{array}
\right),
\end{equation*}
\begin{equation*}
\text{and}\quad A(\bu,\xi)=\sum_{j=1}^n A_j(\bu)\xi_j
= \left(
    \begin{array}[l]{ccc}
      v\cdot \xi & 0& \frac{a_1 }{2}\,^{t}\xi \\
     0 &  v\cdot \xi & \frac{a_2}{2}\,^{t}\xi \\
     2 f' a_1 \, \xi
     &2 f'  a_2\, \xi &  v\cdot \xi I_n
    \end{array}
\right),
\end{equation*}
where $f'$ stands for $f'(|a_1|^2+|a_2|^2)$.
We now have the system:
\begin{equation}
  \label{eq:systhypgeneral}
\partial_t \bu^\eps +\sum_{j=1}^n
  A_j(\bu^\eps)\partial_j \bu^\eps = \frac{\eps}{2} L
  \bu^\eps\quad ;\quad
 \bu^\eps_{\mid t=0}=\bu^\eps_0.
\end{equation}
The matrices $A_j$ are symmetrized by the matrix
\begin{equation*}
  S=\left(
    \begin{array}[l]{cc}
     I_2 & 0\\
     0& \frac{1}{4f'}I_n
    \end{array}
\right),
\end{equation*}
which is symmetric positive if and only if $f'\( |a_1|^2+|a_2|^2\)>0$:
this includes the case of the decofusing cubic nonlinearity
\eqref{eq:gp}, of the cubic-quintic nonlinearity
\eqref{eq:cubicquintic} with $\l>0$, and of the cubic-quintic nonlinearity
\eqref{eq:cubicquintic} with $\l<0$, provided that $|a_1|^2+|a_2|^2 >
|\l|/2$.

\begin{proposition}\label{prop:existencegene}
Assume that $\bu^\eps_0$ is bounded in $X^s$ for all $s>n/2$,
uniformly for $\eps\in [0,1]$, and that there exists $\eps_*>0$
and $\delta>0$ such that
\begin{equation*}
  f'\( |a_0^\eps|^2\)\ge \delta >0,\quad \forall x\in \R^n, \ \forall
  \eps\in [0,\eps_*].
\end{equation*}
Then for $s>n/2+2$, there exist $T_*>0$ and a unique solution
$\bu^\eps \in C([0,T_*];X^s)$  to
\eqref{eq:systhypgeneral} for all $\eps\in [0,\eps_*]$. In addition,
this solution is in $C([0,T_*];X^m)$
for all $m>n/2$, with bounds independent of $\eps\in [0,\eps_*]$.
\end{proposition}
\begin{proof}
Let $s>n/2+2$.
  As usual, the main point consists in obtaining \emph{a priori}
  estimates for the system \eqref{eq:systhypgeneral}, so we shall focus our
  attention on this aspect. We have an \emph{a
  priori} bound for $\bu^\eps$ in $L^\infty$:
\begin{align*}
  \|\bu^\eps(t)\|_{L^\infty}\le &\|\bu^\eps_0\|_{L^\infty}+
  \int_0^t\sum_{j=1}^n \| A_j(\bu^\eps)\d_j \bu^\eps(\tau) \|_{L^\infty}d\tau+
  \int_0^t\|\Delta
  \bu^\eps(\tau)\|_{L^\infty}d\tau\\
\le &\|\bu^\eps_0\|_{L^\infty} + \int_0^t
F\(\|\bu^\eps(\tau)\|_{L^\infty}\)\|\nabla
\bu^\eps(\tau)\|_{H^{s-1}}d\tau\\
& + C\int_0^t \|\Delta
  \bu^\eps(\tau)\|_{H^{s-2}} d\tau.
\end{align*}
We infer:
\begin{equation}\label{eq:bounduLinf0}
 \|\bu^\eps(t)\|_{L^\infty}
\le \|\bu^\eps_0\|_{L^\infty} + \int_0^t G\(\|\bu^\eps(\tau)\|_{X^s}\)
\|\bu^\eps(\tau)\|_{X^s}d\tau.
\end{equation}
To have a closed system of estimates, introduce $P=
(I-\Delta)^{(s-1)/2}\nabla$, so that $\|f\|_{X^s}\thickapprox
\|f\|_{L^\infty} + \|P f\|_{L^2}$. Denote
\begin{equation*}
  \< f,g\>= \int_{\R^n} f(x)\overline{g(x)}dx,
\end{equation*}
the scalar product in $L^2$. Since $S$ is symmetric, we have
\begin{align*}
  \frac{d}{dt}\<SP\bu^\eps(t),P\bu^\eps(t)\>= \<\d_t S
  P\bu^\eps(t),P\bu^\eps(t)\>+2\RE\<S\d_t
  P\bu^\eps(t),P\bu^\eps(t)\>,
\end{align*}
So long as
\begin{equation}
  \label{eq:solong}
  f'\( |a^\eps|^2\)\ge \frac{\delta}{2} >0,
\end{equation}
we have the following set of estimates. First,
\begin{align*}
  \<\d_t S P\bu^\eps(t),P\bu^\eps(t)\>&\le \|\d_t
  S\|_{L^\infty}\|P\bu^\eps(t)\|^2_{L^2}\\
&\le C_\delta \(
  \|\bu^\eps(t)\|_{L^\infty}\)\|\d_t \bu^\eps(t)\|_{L^\infty}
  \|\bu^\eps(t)\|^2_{X^s}.
\end{align*}
Directly from \eqref{eq:systhypgeneral}, we have:
\begin{align*}
  \|\d_t \bu^\eps(t)\|_{L^\infty} &\le C \(
  \|\bu^\eps(t)\|_{L^\infty}\) \|\nabla \bu^\eps(t)\|_{L^\infty} +
  \|\Delta \bu^\eps(t)\|_{L^\infty}\\
 &\le C \(
  \|\bu^\eps(t)\|_{X^s}\) \| \bu^\eps(t)\|_{X^s}.
\end{align*}
Since $SL$ is skew-symmetric, we have
\begin{equation*}
  \RE\< SL P \bu^\eps(t),P\bu^\eps(t)\>=0,
\end{equation*}
which prevents any loss of regularity in the estimates. For the
quasi-linear term involving the matrices $A_j$, we note that since
$SA_j$ is symmetric, commutator estimates (see \cite{DavidJFA}) yield:
\begin{align*}
  \sum_{j=1}^n \< S  P
  \( A_j(\bu^\eps)\partial_j \bu^\eps\),P\bu^\eps(t)\> &\le C\( \|
  \bu^\eps(t) \|_{L^\infty}\) \|P\bu^\eps(t)\|_{L^2}^2
  \|\nabla\bu^\eps(t)\|_{L^\infty}\\
&\le C\( \|
  \bu^\eps(t) \|_{X^s }\)\|P\bu^\eps(t)\|_{L^2}^2.
\end{align*}
Finally, we have:
\begin{equation*}
  \frac{d}{dt} \<SP\bu^\eps(t),P\bu^\eps(t)\>\le C\(
  \|\bu^\eps(t)\|_{X^s}\) \|\bu^\eps(t)\|_{X^s}^2.
\end{equation*}
This estimate, along with \eqref{eq:bounduLinf0}, shows that on a
sufficiently small time interval $[0,T_*]$, with $T_*>0$ independent
of $\eps \in [0,\eps_*]$, \eqref{eq:solong} holds. This
yields the first part of Proposition~\ref{prop:existencegene}.
\smallbreak

The fact that the local existence time does not depend on $s>n/2+2$
follows from the continuation principle based on Moser's calculus and
tame estimates (see e.g. \cite[Section~2.2]{Majda} or
\cite[Section~16.1]{Taylor3}).
\end{proof}
The existence part of Theorem~\ref{theo:cubicquintic} and of
Theorem~\ref{theo:main} when $\phi_{\rm eik}=0$ follows. Indeed, define
$\Phi^\eps$ by
\begin{equation*}
  \Phi^\eps(t)= \phi_0 -\int_0^t \(\frac{1}{2} |v^\eps(\tau)|^2 + f\(
  |a^\eps(\tau)|^2\)\)d\tau.
\end{equation*}
We check that $\d_t ( \nabla \Phi^\eps - v^\eps) =\nabla \d_t
\Phi^\eps - \d_t v^\eps= 0$, so that $\nabla
\Phi^\eps =  v^\eps$, and $(\Phi^\eps,a^\eps)$ solves
\eqref{eq:systemmanuel}.
Finally, uniqueness for \eqref{eq:general} follows from energy
estimates.  If $u^\eps,v^\eps\in C^\infty\cap L^\infty ([0,T_*]\times
\R^n)$ solve \eqref{eq:general}, then $w^\eps := u^\eps
-v^\eps$ satisfies:
\begin{equation*}
  i\eps \d_t w^\eps +\frac{\eps^2}{2}\Delta w^\eps = f\(|u^\eps|^2\)
  u^\eps-f\(|v^\eps|^2\) v^\eps\ ;\ w^\eps_{\mid t=0}=0.
\end{equation*}
We have, for $t\in [0,T_*]$,
\begin{equation*}
  \|w^\eps\|_{L^\infty(0,t;L^2)}\le
  C\(\|u^\eps\|_{L^\infty([0,T_*]\times\R^n)},
\|v^\eps\|_{L^\infty([0,T_*]\times\R^n)}\)
  \|w^\eps\|_{L^1(0,t;L^2)},
\end{equation*}
and Gronwall lemma yields $w^\eps \equiv 0$.

\subsection{Phase-amplitude representation: the case $\phi_{\rm eik}\not =0$}
    
We know consider \eqref{eq:gp}--\eqref{eq:ci} only: the nonlinearity
is exactly cubic. To
take the presence of $V$ and $\phi_{\rm quad}$ into
account, we
proceed as in \cite{CaBKW}: we construct the solution as
$u^\eps = a^\eps e^{i(\phi^\eps +\phi_{\rm eik} )/\eps}$.
The analogue of \eqref{eq:systemmanuel} is:
\begin{equation*}
  \left\{
\begin{aligned}
    \partial_t \Phi^\eps +\frac{1}{2}\left|\nabla
    \Phi^\eps\right|^2 + V+
    |a^\eps|^2-1= 0\quad &; \quad
    \Phi^\eps\big|_{t=0}=\phi_0,\\
\partial_t a^\eps +\nabla \Phi^\eps \cdot \nabla
    a^\eps +\frac{1}{2}a^\eps
\Delta \Phi^\eps  = i\frac{\eps}{2}\Delta
    a^\eps\quad & ;\quad
a^\eps\big|_{t=0}= a^\eps_0\, .
\end{aligned}
\right.
\end{equation*}
Set
$\Phi^\eps= \phi^\eps  +\phi_{\rm eik}$.
The introduction of $\phi_{\rm eik}$ allows us to get rid of the terms
$V_{\rm quad}$ and $\phi_{\rm quad}$, and work in Zhidkov
spaces. The above problem reads, in terms of $(\phi^\eps,a^\eps)$:
\begin{equation}\label{eq:systexact}
\left\{
\begin{aligned}
    \partial_t \phi^\eps +\frac{1}{2}\left|\nabla
    \phi^\eps\right|^2 +\nabla
    \phi_{\rm eik}\cdot \nabla \phi^\eps+V_{\rm lin}+
    |a^\eps|^2-1 &= 0,\\
    \partial_t a^\eps +\nabla \phi^\eps \cdot \nabla
    a^\eps +\nabla \phi_{\rm eik} \cdot \nabla
    a^\eps +\frac{1}{2}a^\eps
\Delta \phi^\eps +\frac{1}{2}a^\eps
\Delta \phi_{\rm eik} &= i\frac{\eps}{2}\Delta a^\eps,\\
\phi^\eps\big|_{t=0}=\phi_{\rm lin}\quad ;\quad a^\eps\big|_{t=0}&=
    a^\eps_0\, .
\end{aligned}
\right.
\end{equation}
Resume the notations of the previous paragraph, with now:
\begin{equation*}
{\bf \Sigma} = \left(
    \begin{array}[l]{c}
       0 \\
       0 \\
       \d_1 V_{\rm lin} \\
      \vdots \\
       \d_n V_{\rm lin}
    \end{array}
\right)
\ , \quad
  \text{and}\quad A(\bu,\xi)=\sum_{j=1}^n A_j(\bu)\xi_j
= \left(
    \begin{array}[l]{ccc}
      v\cdot \xi & 0& \frac{a_1 }{2}\,^{t}\xi \\
     0 &  v\cdot \xi & \frac{a_2}{2}\,^{t}\xi \\
     2  a_1 \, \xi
     &2  a_2\, \xi &  v\cdot \xi I_n
    \end{array}
\right).
\end{equation*}
The system \eqref{eq:systhypgeneral} is replaced by:
\begin{equation}
  \label{eq:systhypV}
\partial_t \bu^\eps +\sum_{j=1}^n
  A_j(\bu^\eps)\partial_j \bu^\eps
  +\nabla \phi_{\rm eik}\cdot \nabla \bu^\eps +
  \tilde M \bu^\eps+{\bf \Sigma}= \frac{\eps}{2} L
  \bu^\eps\quad ;\quad \bu^\eps_{\mid t=0}=\bu^\eps_0,
\end{equation}
where
$\tilde M= \tilde M(t)$ is a matrix depending on time only, since
$\phi_{\rm eik}$ is exactly quadratic in $x$. This aspect seems
necessary in the proof of Proposition~\ref{prop:existence} below. This
explains why we make
Assumptions~\ref{hyp:geom}, and do not content ourselves with general
sub-quadratic potential and initial phase as in \cite{CaBKW}.
The important aspect to notice is that
since the nonlinearity in \eqref{eq:gp} is \emph{exactly} cubic, then
the matrices $A_j$ are symmetrized by a \emph{constant} matrix,
namely:
\begin{equation*}
  S=\left(
    \begin{array}[l]{cc}
     I_2 & 0\\
     0& \frac{1}{4}I_n
    \end{array}
\right).
\end{equation*}
In \cite{CaBKW}, nonlinearities which are cubic \emph{at the origin} were
considered (as in \cite{Grenier98}), and the possibly quadratic phase
$\phi_{\rm eik}$ made the assumption $x a_0^\eps\in L^2(\R^n)$ apparently
necessary, to control the time derivative of the symmetrizer. Of
course, we want to avoid this decay assumption for the Gross--Pitaevskii
equation, so working with a constant symmetrizer  is important.

\begin{proposition}\label{prop:existence}
Assume that $\bu^\eps_0$ is bounded in $X^s$ for all $s>n/2$,
uniformly for $\eps\in [0,1]$. Then for $s>n/2+2$, there exist $T_*\in ]0,T]$,
independent of $\eps\in [0,1]$, and a unique solution
$\bu^\eps \in C([0,T_*];X^s)$  to
\eqref{eq:systhypV}. In addition, this solution is in $C([0,T_*];X^m)$
for all $m>n/2$, with bounds independent of $\eps\in [0,1]$.
\end{proposition}
\begin{proof}[Sketch of the proof]
The proof follows the same lines as the proof of
Proposition~\ref{prop:existencegene}, so we shall only point out the
differences.

Let $s>n/2+2$. By construction, the
  operator $\d_t +\nabla \phi_{\rm eik}\cdot \nabla$ is a transport
  operator along the characteristics associated to $\phi_{\rm eik}$,
  which do not intersect for $t\in [0,T]$. Therefore, we have an \emph{a
  priori} bound for $\bu^\eps$ in $L^\infty$:
\begin{align}
  \|\bu^\eps(t)\|_{L^\infty}\le &\|\bu^\eps_0\|_{L^\infty}+
  \int_0^t\sum_{j=1}^n \| A_j(u)\d_j u(\tau) \|_{L^\infty}d\tau\notag\\
& +  \int_0^t\Big( C
  \|\bu^\eps(\tau)\|_{L^\infty} +\|{\bf \Sigma}(\tau)\|_{L^\infty}+  \|\Delta
  \bu^\eps(\tau)\|_{L^\infty}\Big) d\tau\notag \\
\le & \|\bu^\eps_0\|_{L^\infty}+ C\int_0^t\( 1+
  \|\bu^\eps(\tau)\|_{X^s}\) \|\bu^\eps(\tau)\|_{X^s} d\tau + \|{\bf
  \Sigma}\|_{L^\infty([0,T];X^s)}.  \label{eq:bounduLinf}
\end{align}
To have a closed system of estimates, resume the operator $P=
(I-\Delta)^{(s-1)/2}\nabla$, so that $\|f\|_{X^s}\thickapprox
\|f\|_{L^\infty} + \|P f\|_{L^2}$. We have
\begin{align*}
  \frac{d}{dt}\<SP\bu^\eps(t),P\bu^\eps(t)\>= 2\RE\<S\d_t
  P\bu^\eps(t),P\bu^\eps(t)\>,
\end{align*}
since $S$ is constant symmetric. Since $SL$ is skew-symmetric, we have
\begin{equation*}
  \RE\< SL P \bu^\eps(t),P\bu^\eps(t)\>=0.
\end{equation*}
For the
quasi-linear term involving the matrices $A_j$, we note that since
$SA_j$ is symmetric, commutator estimates yield:
\begin{align*}
  \sum_{j=1}^n \< S  P
  \( A_j(\bu^\eps)\partial_j \bu^\eps\),P\bu^\eps(t)\> \le C\( \|
  \bu^\eps(t) \|_{X^s }\)\|P\bu^\eps(t)\|_{L^2}^2.
\end{align*}
Next, write
\begin{align*}
  \< S P \(\nabla \phi_{\rm eik}\cdot \nabla  \bu^\eps(t)\),
  P\bu^\eps(t)\> =& \< S  \nabla \phi_{\rm eik}\cdot \nabla P \bu^\eps(t),
  P\bu^\eps(t)\>\\
& + \< S  [P,\nabla \phi_{\rm eik}\cdot \nabla ] \bu^\eps(t),
  P\bu^\eps(t)\>.
\end{align*}
The first term of the right-hand side is estimated thanks to an
integration by parts:
\begin{align*}
  2\RE\< S  \nabla \phi_{\rm eik}\cdot \nabla P \bu^\eps(t),
  P\bu^\eps(t)\> &= \int S  \nabla \phi_{\rm eik}(t,x)\cdot \nabla |P
  \bu^\eps(t,x)|^2dx\\
& = - \int S  \Delta \phi_{\rm eik}(t,x) |P \bu^\eps(t,x)|^2dx.
\end{align*}
For the second term, we notice that $[P,\nabla \phi_{\rm eik}\cdot
\nabla ] = \psi\nabla$, where $\psi=\psi(t,D)$ is a
pseudo-differential operator in $x$, of order $s-1$, depending
smoothly of $t\in [0,T]$. Therefore,
\begin{align*}
  2\RE\< S P \(\nabla \phi_{\rm eik}\cdot \nabla  \bu^\eps(t)\),
  P\bu^\eps(t)\> \lesssim \|\bu^\eps(t)\|_{X^s}^2.
\end{align*}
The fact that $\tilde M$ is independent of $x$ is crucial here, to
ensure that $P(\tilde M \bu^\eps)\in L^2$ for $\bu^\eps\in X^s$. If
$\tilde M$ depended on $x$, that is if $\phi_{\rm eik}$ was not a
polynomial of order at most two, the low frequencies might be a
problem at this step of the proof. Finally, we have:
\begin{equation*}
  \frac{d}{dt} \<SP\bu^\eps(t),P\bu^\eps(t)\>\le C\(
  \|\bu^\eps(t)\|_{X^s}\) \|\|\bu^\eps(t)\|_{X^s}^2.
\end{equation*}
This estimate, along with \eqref{eq:bounduLinf}, yields the first part
of Proposition~\ref{prop:existence}. We conclude like in the proof of
Proposition~\ref{prop:existencegene}.
\end{proof}
The existence part of Theorem~\ref{theo:main} follows from the above
result, by setting
\begin{equation*}
  \phi^\eps(t) = \phi_{\rm lin} -\int_0^t \(
    \frac{1}{2}\left|v^\eps(\tau)\right|^2 +\nabla
    \phi_{\rm eik}(\tau)\cdot v^\eps(\tau)+V_{\rm lin}(\tau)+
    |a^\eps(\tau)|^2-1\)d\tau.
\end{equation*}
Finally, uniqueness for \eqref{eq:gp}--\eqref{eq:ci} follows from energy
estimates.  If $u^\eps,v^\eps\in C^\infty\cap L^\infty ([0,T_*]\times
\R^n)$ solve \eqref{eq:gp}--\eqref{eq:ci}, then $w^\eps := u^\eps
-v^\eps$ satisfies:
\begin{equation*}
  i\eps \d_t w^\eps +\frac{\eps^2}{2}\Delta w^\eps = (V-1)w^\eps
  +|u^\eps|^2 u^\eps-|v^\eps|^2 v^\eps\ ;\ w^\eps_{\mid t=0}=0.
\end{equation*}
We have, for $t\in [0,T_*]$,
\begin{equation*}
  \|w^\eps\|_{L^\infty(0,t;L^2)}\lesssim \(
  \|u^\eps\|_{L^\infty([0,T_*]\times\R^n)}^2 +
  \|v^\eps\|_{L^\infty([0,T_*]\times\R^n)}^2\)\|w^\eps\|_{L^1(0,t;L^2)},
\end{equation*}
and Gronwall lemma yields $w^\eps \equiv 0$.

\subsection{On the Hamiltonian structure}
\label{sec:hamilton}

When  $V=V(x)$ is time-independent, \eqref{eq:gp} formally has a
Hamiltonian structure, with
\begin{equation*}
  H= \frac{1}{2}\|\eps \nabla u^\eps(t)\|_{L^2}^2
  +\int_{\R^n}V(x)|u^\eps(t,x)|^2dx +
  \frac{1}{2}\left\||u^\eps(t)|^2-1\right\|_{L^2}^2 .
\end{equation*}
When $V\equiv 0$, this structure is used in
\cite{PG05} to prove the global existence of solutions in the energy
space. On the other hand, suppose that $V$ is, say, harmonic:
\begin{equation*}
  V(x)=\sum_{j=1}^n \lambda_j x_j^2,
\end{equation*}
where the constants $\lambda_j\ge 0$ are not all equal to zero. Then
necessarily, $H$ is infinite: suppose for instance that
$\lambda_1>0$. Then if  $\d_{x_1}u^\eps(t,\cdot),x_1u^\eps (t,\cdot)
\in L^2(\R^n)$, the 
uncertainty principle (a simple integration by parts, plus
Cauchy--Schwarz inequality in this case) yields:
\begin{equation*}
  u^\eps (t,\cdot) \in L^2(\R^n).
\end{equation*}
Therefore,  the constraint
 $|u^\eps(t,\cdot)|^2 -1\in
L^2(\R^n)$ cannot be satisfied, for otherwise,
 $1=1-|u^\eps(t,\cdot)|^2 +|u^\eps(t,\cdot)|^2 \in L^2(\R^n)+
 L^1(\R^n)$.
\smallbreak

Similarly, assume that $V\equiv 0$, but $\phi_{\rm quad}\not =0$:
rapid quadratic oscillations are present in the initial data. We have
\begin{equation*}
  \eps \nabla u^\eps_{\mid t=0} = \(\eps \nabla a_0^\eps +
  i a_0^\eps \nabla \phi_0\)e^{i\phi_0/\eps}.
\end{equation*}
Therefore, the above quantity is in $L^2$ provided that $\nabla
a_0^\eps ,a_0^\eps \nabla \phi_{\rm quad}\in
L^2(\R^n)$. If for instance $\phi_{\rm quad}(x) = cx_1^2$ with
$c\not =0$, the last assumption means that $x_1a_0^\eps \in
L^2(\R^n)$, which brings us back to the previous discussion.
\smallbreak

We shall see in Section~\ref{sec:hamil} that if $\phi_{\rm
  eik}\not \equiv 0$,
and if $a_0^\eps\in X^\infty$ is such that
\begin{equation*}
\left|a_0^\eps \right|^2-1\in L^2(\R^n),
\end{equation*}
then the last constraint present in $H$ is not propagated in
general. In small time at least, one has generically
\begin{equation*}
\left|u^\eps(t,\cdot) \right|^2-1\not\in L^2(\R^n).
\end{equation*}

\section{Semi-classical analysis}
\label{sec:semi}
We now complete the proof of Theorem~\ref{theo:main}. The end of the
proof of Theorem~\ref{theo:cubicquintic} follows essentially the same
lines, so we omit it. The main adaptation is due to the fact that when
the nonlinearity is not exactly cubic, the symmetrizer $S$ is not
constant. We refer to \cite{Grenier98} or \cite{CaBKW}, to see that
the proof below is easily adapted.
\smallbreak

Introduce $(\phi,a)$, solution to \eqref{eq:systexact} with $\eps=0$,
that is
\begin{equation}\label{eq:systexact0}
\left\{
\begin{aligned}
    \partial_t \phi +\frac{1}{2}\left|\nabla
    \phi\right|^2 +\nabla
    \phi_{\rm eik}\cdot \nabla \phi+V_{\rm lin}+
    |a|^2-1 &= 0\ ;\ \phi\big|_{t=0}=\phi_{\rm lin}\, ,\\
    \partial_t a +\nabla \phi \cdot \nabla
    a+\nabla \phi_{\rm eik} \cdot \nabla
    a +\frac{1}{2}a
\Delta \phi +\frac{1}{2}a
\Delta \phi_{\rm eik} &= 0\ ; \ a\big|_{t=0}= a_0\, .
\end{aligned}
\right.
\end{equation}
It is a particular case of Proposition~\ref{prop:existence} that
\eqref{eq:systexact0} has a unique solution, such that $a,\nabla \phi \in
C([0,T_*];X^s)$ for all $s>n/2$.

\begin{proposition}\label{prop:error1}
  Under the assumptions of Theorem~\ref{theo:main}, let
  $(\phi^\eps,a^\eps)$ and $(\phi,a)$ be given by \eqref{eq:systexact}
  and \eqref{eq:systexact0} respectively. For all $s>n/2$, there
  exists $C_s$ such that
  \begin{equation*}
    \|\nabla(\phi^\eps-\phi)\|_{L^\infty([0,T_*];X^s)} +
    \|a^\eps-a\|_{L^\infty([0,T_*];X^s)} \le C_s\eps.
  \end{equation*}
\end{proposition}
\begin{proof}[Sketch of the proof]
  We shall give the outline of the proof, since it is very similar to
  the case of Sobolev spaces \cite{CaBKW}. The differences are those
  pointed out in the proof of
  Proposition~\ref{prop:existence}. Resuming the notations of
  \S\ref{sec:general}, set
  \begin{equation*}
  \bu = \left(
    \begin{array}[l]{c}
       \RE a \\
       \IM a \\
       \d_1 \phi \\
      \vdots \\
       \d_n \phi
    \end{array}
\right),\quad
\bw_0^\eps  = \left(
    \begin{array}[l]{c}
       \RE (a_0^\eps -a_0) \\
       \IM (a_0^\eps -a_0)\\
       0 \\
      \vdots \\
       0
    \end{array}
\right)
\end{equation*}
Denoting $\bw^\eps = \bu^\eps -\bu$, \eqref{eq:systhypV} yields:
 \begin{equation*}
\left\{
  \begin{aligned}
    &\partial_t \bw^\eps +\sum_{j=1}^n
  \(A_j(\bu^\eps)\partial_j \bu^\eps - A_j(\bu)\partial_j \bu\)
  +\nabla \phi_{\rm eik}\cdot \nabla \bw^\eps +
  \tilde M \bw^\eps= \frac{\eps}{2} L
  \bw^\eps + \frac{\eps}{2} L
  \bu\, ,\\
 &\bw^\eps_{\mid t=0}=\bw^\eps_0.
  \end{aligned}
\right.
\end{equation*}
We know by Proposition~\ref{prop:existence} that $\bu^\eps$ and $\bu$
are bounded in $C([0,T_*];X^s)$ for all $s>n/2$. The source term ${\bf
  \Sigma}$ in \eqref{eq:systhypV} is now replaced by $\frac{\eps}{2} L
  \bu$, which is of order $\O(\eps)$ in $C([0,T_*];X^s)$, and we have
  easily, for $s>n/2$ and $t\in [0,T_*]$:
  \begin{equation*}
    \| \bw^\eps \|_{L^\infty([0,t];X^s)} \le \|\bw^\eps_0\|_{X^s} +\O(\eps)+
    \int_0^t \| \bw^\eps (\tau)\|_{X^s}d\tau.
  \end{equation*}
The proposition follows from Gronwall lemma.
\end{proof}
\begin{remark}\label{rem:global}
  Note that for the time $T_*$ in Proposition~\ref{prop:error1} (as well
  as in Proposition~\ref{prop:error2} below), we can pick the
  life-span of $(\phi,a)$, the solution of
  \eqref{eq:systexact0}. Indeed, the error estimate and the standard
  continuity argument  show that $(\phi^\eps,a^\eps)$ cannot blow-up
  as long as $(\phi,a)$ remains smooth, provided that $\eps$ is chosen
  sufficiently small. In particular, if $(\phi,a)$ remains smooth
  globally in time, then for any $\tau>0$, we can find $\eps(\tau)>0$
  such that Proposition~\ref{prop:error1} and
  Proposition~\ref{prop:error2} below remain valid on $[0,\tau]$ for
  $\eps \in ]0,\eps(\tau)]$. On the other hand, one must not expect
  $T_*=\infty$ in general: the solution to
\eqref{eq:euler} may not remain smooth for all time. See
\cite{Sideris}. 

\end{remark}
\begin{corollary}\label{cor:linf}
  There
  exists $C$ such that for all $t\in [0,T_*]$,
  \begin{equation*}
    \|\phi^\eps(t,\cdot)-\phi (t,\cdot)\|_{L^\infty}  \le C\eps t.
  \end{equation*}
\end{corollary}
\begin{proof}
  Set $w_\phi^\eps = \phi^\eps-\phi $. It satisfies
  \begin{equation*}
    (\d_t + \nabla \phi_{\rm eik}\cdot \nabla)w_\phi^\eps =
    \frac{1}{2}\( |\nabla \phi|^2 - |\nabla \phi^\eps|^2\) +
    |a|^2-|a^\eps|^2\ ;\ w_{\phi \mid t=0}^\eps =0.
  \end{equation*}
By Proposition~\ref{prop:error1}, the right hand side is
$\O(\eps)$ in $L^\infty$. Integration along
the characteristics associated to $\d_t + \nabla \phi_{\rm eik}\cdot
\nabla$ (see Remark~\ref{rema:transport}) yields the result.
\end{proof}
The first estimate \eqref{eq:BKWVtpetit} of Theorem~\ref{theo:main}
follows easily:
\begin{align*}
  u^\eps - ae^{i\phi/\eps} &= a^\eps e^{i\phi^\eps/\eps} -
  ae^{i\phi/\eps} = \( a^\eps -a\) e^{i\phi^\eps/\eps} +
  a\(e^{i\phi^\eps/\eps} - e^{i\phi/\eps}\)\\
&= \O(\eps) +  a e^{i(\phi^\eps +\phi)/(2\eps)}2i\sin \(
  \frac{\phi^\eps -\phi}{2\eps}\) = \O(\eps) + \O (t),
\end{align*}
where the $\O(\cdot)$'s stand for estimates in $L^\infty([0,T_*]\times
\R^n)$.
\smallbreak

To improve \eqref{eq:BKWVtpetit} to \eqref{eq:BKWVtgrand}, we need the
next term in the asymptotic expansion of $(\phi^\eps,a^\eps)$ in terms
of powers of $\eps$. Introduce the system:
\begin{equation}
  \label{eq:systlinear}
\left\{
  \begin{aligned}
  \d_t \phi^{(1)}+ \nabla
    (\phi_{\rm eik}+\phi)\cdot \nabla \phi^{(1)}
    +2\RE \(\overline{a}a^{(1)}\) &= 0\ ;\ \phi^{(1)}_{\mid t=0}=0.\\
  \d_t a^{(1)}+ \nabla
    (\phi_{\rm eik}+\phi)\cdot \nabla  a^{(1)}+\nabla \phi^{(1)}\cdot
    \nabla a + \frac{1}{2}&a^{(1)}\Delta (\phi_{\rm
    eik}+\phi)\\
+\frac{1}{2}a\Delta  \phi^{(1)}&=\frac{i}{2}\Delta a \ ;\
    a^{(1)}_{\mid t=0}=a_1.
  \end{aligned}
\right.
\end{equation}
It is easy to see that 
this linear system has a unique classical solution such that  
$a^{(1)},\nabla \phi^{(1)}\in C([0,T_*];X^s)$ for all
$s>n/2$. Reasoning as in the proof of Corollary~\ref{cor:linf}, we see
that we have also $\phi^{(1)}\in C([0,T_*];X^s)$. Moreover, mimicking
the proofs of Proposition~\ref{prop:error1} and 
Corollary~\ref{cor:linf}, we have the following result, whose proof is
left out:
\begin{proposition}\label{prop:error2}
  Let
  $(\phi^\eps,a^\eps)$, $(\phi,a)$ and $(\phi^{(1)},a^{(1)})$ be given
  by \eqref{eq:systexact},
 \eqref{eq:systexact0} and \eqref{eq:systlinear} respectively. Denote
  $r_0^\eps = a_0^\eps -a_0 -\eps a_1$.  For all
  $s>n/2+2$,
  \begin{equation*}
    \|\nabla(\phi^\eps-\phi-\eps \phi^{(1)})\|_{L^\infty([0,T_*];X^s)} +
    \|a^\eps-a-\eps a^{(1)} \|_{L^\infty([0,T_*];X^s)} \le \widetilde
C_s\(\eps^2    + \|r_0^\eps\|_{X^s}\).
  \end{equation*}
In addition, there exists $\widetilde C$ such that if $s>n/2+2$,
\begin{equation*}
  \|\phi^\eps-\phi-\eps \phi^{(1)}\|_{L^\infty([0,T_*]\times \R^n)}\le
  \widetilde C \(\eps^2    + \|r_0^\eps\|_{X^s}\).
\end{equation*}
\end{proposition}
We can now complete the proof of Theorem~\ref{theo:main}:
\begin{align*}
  u^\eps - ae^{i\phi^{(1)}}e^{i\phi/\eps} &= a^\eps e^{i\phi^\eps/\eps} -
  ae^{i(\phi+\eps \phi^{(1)})/\eps}\\
& = \( a^\eps -a\) e^{i\phi^\eps/\eps} +
  a\(e^{i\phi^\eps/\eps} - e^{i(\phi+\eps \phi^{(1)})/\eps}\)\\
&= \O(\eps) +  a e^{i(\phi^\eps +\phi+\eps \phi^{(1)})/(2\eps)}2i\sin \(
  \frac{\phi^\eps -\phi- \eps \phi^{(1)}}{2\eps}\)\\
& = \O(\eps) + \O
  \(\frac{\|r_0^\eps\|_{X^s}}{\eps}\).
\end{align*}
This yields \eqref{eq:BKWVtgrand}, along with Remark~\ref{rema:mieux}.
\begin{remark}\label{rem:arb}
  Following the same lines, we see that if $a_0^\eps$ is known up to order
  $\O(\eps^{N+1})$ in $X^s$ for some $s>n/2+2$, $N\in \N$, then we can
  construct
  an approximate solution $v_N^\eps$ such that
  \begin{equation*}
    \| u^\eps -v_N^\eps\|_{L^\infty([0,T_*];X^s)}=\O\(\eps^N\).
  \end{equation*}
\end{remark}
To conclude this paragraph, we note that if we know that the initial
corrector $a_1$ is not only in $X^\infty$, but in $H^\infty$, then
Theorem~\ref{theo:main} becomes more precise.
\begin{corollary}\label{cor:hs}
Under the same assumptions as in Theorem~\ref{theo:main},
suppose moreover that $a_1\in H^\infty$, and
\begin{equation*}
    \left\|a_0^\eps -a_0 -\eps a_1\right\|_{H^s}
    =\O\(\delta^\eps\),\quad \forall s>0,\quad  \text{with
    }\delta^\eps =o(\eps) .
  \end{equation*}
Then \eqref{eq:BKWVtgrand} can be improved to:
\begin{equation}\label{eq:BKWVtgrandHS}
    \sup_{t\in [0,T_*]}\left\|
    u^\eps(t,\cdot) - a(t,\cdot)e^{i\phi^{(1)}(t,\cdot)}
    e^{i(\phi(t,\cdot)+\phi_{\rm
    eik}(t,\cdot))/\eps}\right\|_{L^\infty\cap L^2}=\O\(\eps +
    \frac{\delta^\eps}{\eps}\) .
  \end{equation}
\end{corollary}
Essentially, one just has to notice that the error estimates in
Propositions~\ref{prop:error1} and \ref{prop:error2} can then be
measured in $H^s$ instead of $X^s$. Note also that in
\eqref{eq:BKWVtgrandHS}, it may happen that none of the two functions
is in $L^2$.

\section{Examples when $\phi_{\rm eik}\equiv 0$}
\label{sec:modone}

In this paragraph, we consider \eqref{eq:gp}--\eqref{eq:ci}, and we
assume $\phi_{\rm eik}\equiv 0$.

\subsection{An example from \cite{ColinSoyeur}}
\label{sec:colinsoyeur}

As an application, we can recover and improve the result of
\cite{ColinSoyeur}, in the case of the whole space (the space variable
$x$ lies in a bounded domain in \cite{ColinSoyeur}). Assume  that
\begin{equation*}
  a_0^\eps(x) = a_0(x)=e^{i\theta_0(x)},\quad \theta_0\in
  H^\infty(\R^n;\R)\quad
  ;\quad \phi_0=V=0.
\end{equation*}
That is, we consider:
\begin{equation*}
  i\eps \d_t u^\eps +\frac{\eps^2}{2}\Delta u^\eps =\(
  |u^\eps|^2-1\)u^\eps\quad
  ;\quad
u^\eps(0,x)=e^{i\theta_0(x)}.
\end{equation*}
Then $a_0^\eps = a_0\in X^\infty$, and we see that:
\begin{itemize}
\item $\phi\equiv 0$ and $a$ is independent of time:
  $a(t,x)=a_0(x)=e^{i\theta_0(x)}$.
\item $\phi^{(1)}$ solves
  \begin{equation*}
    \d_t^2 \phi^{(1)} = \IM\( \overline a \Delta a\),
  \end{equation*}
so that $\theta(t,x) \defn \phi^{(1)}(t,x)+\theta_0(x)$ solves:
\begin{equation*}
  \(\d_t^2 -\Delta \)\theta =0\quad ;\quad
  \theta(0,x)=\theta_0(x)\quad;\quad \d_t \theta(0,x)=0.
\end{equation*}
\end{itemize}
Note that $(\phi,a)$ remains smooth for all time, so we can take
$T_*$ arbitrarily large (see Remark~\ref{rem:global}).
Since from Theorem~\ref{theo:main} and the above corollary,
\begin{equation*}
  \sup_{t\in [0,T_*]}\|u^\eps(t,\cdot) - a(t,\cdot) e^{i
    \phi^{(1)}(t,\cdot)}\|_{L^\infty\cap L^2}=
\sup_{t\in [0,T_*]}\|u^\eps(t,\cdot) -  e^{i
    \theta(t,\cdot)}\|_{L^\infty\cap L^2} =\O(\eps),
\end{equation*}
where $T_*>0$ is arbitrary. We recover
\cite[Theorem~2]{ColinSoyeur} in the case of the whole space, with no
restriction on the space dimension, and a precised error
estimate. Note also that in view of
Remark~\ref{rem:arb}, we can justify \cite[Proposition~5]{ColinSoyeur}
(giving the $\eps$-order corrector for $u^\eps$), and get a complete
asymptotic expansion for $u^\eps$.

\subsection{When $|a_0^\eps|^2-1\in L^2$}\label{sec:gpi}
As in Corollary~\ref{cor:hs}, assume that \eqref{eq:DACI} is precised to
\begin{equation*}
    \left\|a_0^\eps -a_0 -\eps a_1\right\|_{H^s}
    =o(\eps),\quad \forall s>0.
  \end{equation*}
where $a_0\in X^\infty$ and $a_1\in H^\infty$. Assume moreover that
\begin{equation*}
  |a_0|^2-1\in L^2(\R^n).
\end{equation*}
Then \eqref{eq:systemmanuel} yields:
\begin{align*}
  \frac{d}{dt}&\left\||a^\eps(t)|^2-1\right\|_{L^2}^2 = 4\int_{\R^n}
  \left||a^\eps(t,x)|^2-1\right|
\RE \(\overline a^\eps(t,x) \d_t a^\eps(t,x)\) dx \\
& \lesssim
  \left\||a^\eps(t)|^2-1\right\|_{L^2}\|a^\eps\|_{L^\infty}
 \( \|\nabla \Phi^\eps\cdot \nabla a^\eps\|_{L^2} + \|
  a^\eps \Delta \Phi^\eps\|_{L^2} +\|
  \Delta a^\eps \|_{L^2}  \)\\
& \lesssim
  \left\||a^\eps(t)|^2-1\right\|_{L^2}\|a^\eps\|_{X^s}^2\(\|\nabla
  \Phi^\eps\|_{L^\infty}+ \|\Delta \Phi^\eps\|_{L^2}+1\),
\end{align*}
where we consider $s>n/2+2$. Therefore,
Proposition~\ref{prop:existencegene} shows that
\begin{equation*}
 |u^\eps|^2-1\in C([0,T_*];L^2(\R^n)).
\end{equation*}
Note that this property holds even if $V=V_{\rm lin}\not =0$.

\subsection{When $a_0^\eps(x)\sim 1$ as $|x|\to \infty$}
In a spirit similar to \cite{LinZhang} (where the authors choose
$\theta_0\equiv
0$), assume that $V=0$, $\phi_0(x)= v^\infty\cdot x$ for
some $v^\infty\in \R^n$, and
\begin{equation*}
    \left\|a_0^\eps - e^{i\theta_0(x)} -\eps a_1\right\|_{H^s}
    =\O\(\delta^\eps\),\quad \forall s>0,\quad \text{where }\theta_0\in
    H^\infty\text{ is real-valued}.
\end{equation*}
Then as in \S\ref{sec:colinsoyeur}, we compute:
\begin{equation*}
  \phi(t,x)= v^\infty\cdot x -\frac{|v^\infty|^2}{2}t\quad ;\quad
  a(t,x)= a_0\(x-v^\infty t\)=e^{i\theta_0\(x-v^\infty t\)}.
\end{equation*}
We also note that $T_*>0$ can be taken arbitrarily large. In addition,
we check that $\phi^{(1)}$ is such that $\widetilde \phi^{(1)}
(t,y)=\phi^{(1)} (t,x+v^\infty t)$ solves:
\begin{equation*}
  \(\d_t^2-\Delta\) \widetilde \phi^{(1)} = \IM \( \overline a_0\Delta
  a_0\) = \Delta \theta_0\ \  ;\ \  \widetilde \phi^{(1)} (0,x)= 0
  \ ;\ \d_t\widetilde \phi^{(1)} (0,x)= -2\RE \(\overline a_0
  a_1\).
\end{equation*}
Therefore, Corollary~\ref{cor:hs} yields
\begin{equation*}
  \sup_{t\in [0,T]} \left\| u^\eps(t,\cdot) -
  e^{i\theta(t,\cdot)}e^{i\phi(t,\cdot)/\eps}\right\|_{L^\infty\cap L^2}=
  \O\(\eps+\frac{\delta^\eps}{\eps} \),
\end{equation*}
where $\theta$ is given by $\theta(t,x)= \widetilde
\theta(t,y)\big|_{y= x-v^\infty t}$, with:
\begin{equation*}
  \(\d_t^2-\Delta\) \widetilde\theta =0\quad ;\quad \widetilde\theta_{\mid
  t=0}=\theta_0\ ;\ \d_t  \widetilde\theta_{\mid
  t=0}=-2\RE \(\overline a_0 a_1\).
\end{equation*}

\section{Time propagation of the condition at infinity: $\phi_{\rm
eik}\not =0$}
\label{sec:hamil}

In this section, we assume that $|a_0^\eps|^2-1\in L^2(\R^n)$, and aim
at understanding how this condition is propagated on the time interval
$[0,T_*]$ when $\phi_{\rm
eik}\not =0$. Essentially, we have $|u^\eps(t,\cdot)|^2-1\in  L^2(\R^n)$
for $t\in [0,T_*]$ if and only if $\phi_{\rm eik}\equiv 0$. The
function $\phi_{\rm eik}$ is identically zero if and only if $V_{\rm
  quad}=\phi_{\rm quad}=0$: that case was developed in
\S\ref{sec:gpi}. We compute
\begin{align*}
  \frac{d}{dt}\left\| |u^\eps(t)|^2-1\right\|_{L^2}^2 \le 4 \left\|
  |u^\eps(t)|^2-1\right\|_{L^2}
  \|a^\eps(t)\|_{L^\infty}\|\d_t a^\eps(t)\|_{L^2}.
\end{align*}
In the above estimate, we assumed that $\d_t a^\eps(t,\cdot)\in
L^2$. Let us now examine this condition. In view of
Proposition~\ref{prop:existence}, we know
that all the terms in the second equation of \eqref{eq:systexact} are
in $L^2(\R^n)$, except
possibly $\d_t a^\eps$, $\nabla \phi_{\rm eik}\cdot \nabla a^\eps$ and
$a^\eps\Delta \phi_{\rm eik}$. Therefore if $\phi_{\rm eik}\equiv 0$,
we infer that $|u^\eps(t,\cdot)|^2-1 \in L^2(\R^n)$ for all $t\in
[0,T_*]$.
\smallbreak

Assume now that $\phi_{\rm eik}$ is not zero.
To gather the terms $\d_t a^\eps$ and $\nabla \phi_{\rm eik}\cdot
\nabla a^\eps$ together, consider the change of variable of
Remark~\ref{rema:transport}, and set
\begin{equation*}
  \widetilde a^\eps(t,y) =a^\eps (t,x(t,y)).
\end{equation*}
Since the Jacobi determinant $\operatorname{det} \nabla_y x(t,y)>0$ is
bounded from above, and from below away from zero for $t\in
[0,T_*]\subset [0,T]$,  $\d_t a^\eps(t,\cdot)$ and $\d_t \widetilde
a^\eps(t,\cdot)$ are simultaneously in $L^2(\R^n)$. Given $\Delta
  \phi_{\rm eik}$ is a function of time only, we have
\begin{equation*}
  \d_t \widetilde a^\eps = -\frac{1}{2}\widetilde a^\eps \Delta
  \phi_{\rm eik} +C([0,T_*];L^2).
\end{equation*}
We are in a case where $\widetilde a^\eps \Delta
  \phi_{\rm eik} \not \in L^2$. To overcome this issue, consider
  \begin{equation*}
     \left\| \left|u^\eps(t) e^{\frac{1}{2} \int_0^t \Delta
  \phi_{\rm eik} (\tau)d\tau}\right|^2-1\right\|_{L^2}^2 =
\left\| \left|a^\eps(t) e^{\int_0^t
  \operatorname{Tr}Q (\tau)d\tau}\right|^2-1\right\|_{L^2}^2,
  \end{equation*}
where $Q$ is given by Lemma~\ref{lem:hj}. For $t\in [0,T_*]$, this
quantity is equivalent to:
\begin{equation*}
\left\| \left|\widetilde a^\eps(t) e^{\int_0^t
  \operatorname{Tr}Q (\tau)d\tau}\right|^2-1\right\|_{L^2}^2.
  \end{equation*}
We have:
\begin{align*}
  \frac{d}{dt}  \left\| \left|\widetilde a^\eps(t) e^{\int_0^t
  \operatorname{Tr}Q (\tau)d\tau}\right|^2-1\right\|_{L^2}^2\le
C& \left\| \left|\widetilde a^\eps(t) e^{\int_0^t
  \operatorname{Tr}Q (\tau)d\tau}\right|^2-1\right\|_{L^2}\|\widetilde
  a^\eps(t)\|_{L^\infty}\times\\
&\times \left\| \d_t \widetilde a^\eps(t) +
  \frac{1}{2}\widetilde a^\eps(t) \Delta \phi_{\rm eik}(t)\right\|_{L^2}.
\end{align*}
We infer that $\left|\widetilde a^\eps(t) e^{\int_0^t
  \operatorname{Tr}Q (\tau)d\tau}\right|^2-1 \in C([0,T_*];L^2)$,
  hence
$$\left|u^\eps(t) e^{\int_0^t
  \operatorname{Tr}Q (\tau)d\tau}\right|^2-1 \in
  C([0,T_*];L^2).$$
Morally, for $t\in [0,T_*]$, the modulus of $u^\eps$ goes to $\exp(-\int_0^t
  \operatorname{Tr}Q (\tau)d\tau)$ as $|x|\to \infty$. We conclude by
  some examples that illustrate this analysis.
\smallbreak

\noindent\emph{Example 1}. Consider the case where $\phi_{\rm quad}=0$, and
  $V_{\rm quad}(x) = \omega^2 \frac{|x|^2}{2}$ is an isotropic
  harmonic potential ($\omega>0$). Then we compute
  \begin{align*}
   & \phi_{\rm eik}(t,x) = -\omega\frac{|x|^2}{2}\tan (\omega t), \quad
    t\in [0,T]\subset \left[ 0,\frac{\pi}{2\omega}\right[, \\
\text{and}\qquad &  \exp\(-\int_0^t
  \operatorname{Tr}Q (\tau)d\tau\) =\exp\(\frac{n\omega}{2}\int_0^t
  \tan (\omega \tau)d\tau\) = \( \cos (\omega t)\)^{-n/2}.
\end{align*}
Therefore, the ``limit of the modulus of $u^\eps$ at infinity'' grows
at time evolves. If in Proposition~\ref{prop:existence}, we can take
$T_*$ arbitrarily close to $\pi/(2\omega)$, this suggests that there
is some sort of ``blow-up at infinity'' at $t$ approaches
$\pi/(2\omega)$.
\smallbreak

\noindent\emph{Example 2}. Consider the case where $\phi_{\rm quad}=0$, and
  $V_{\rm quad}(x) = -\omega^2 \frac{|x|^2}{2}$ is an isotropic
  \emph{repulsive} harmonic potential ($\omega>0$). We have
  \begin{align*}
    &\phi_{\rm eik}(t,x) = \omega\frac{|x|^2}{2}\tanh (\omega t), \quad
    t\in [0,+\infty[,\\
 \text{and}\qquad &   \exp\(-\int_0^t
\operatorname{Tr}Q (\tau)d\tau\) =\( \cosh (\omega t)\)^{-n/2}.
\end{align*}
Therefore, the ``limit of the modulus of $u^\eps$ at infinity'' decays
at time evolves.
\smallbreak

\noindent\emph{Example 3}. Consider the case  $\phi_{\rm
  quad}=-|x|^2/2$, and
  $V_{\rm quad}(x) = 0$. We compute
  \begin{equation*}
    \phi_{\rm eik}(t,x) = \frac{|x|^2}{2(t-1)}, \quad
    t\in [0,1[, \quad \text{and }
  \exp\(-\int_0^t
  \operatorname{Tr}Q (\tau)d\tau\) =\( 1-t\)^{-n/2}.
\end{equation*}
This case is similar to the first example.
\smallbreak

\noindent\emph{Example 4}. Consider the case $\phi_{\rm
  quad}=|x|^2/2$, and
  $V_{\rm quad}(x) = 0$. We have
  \begin{equation*}
    \phi_{\rm eik}(t,x) = \frac{|x|^2}{2(t+1)}, \quad
    t\in [0,+\infty[, \quad \text{and }
  \exp\(-\int_0^t
  \operatorname{Tr}Q (\tau)d\tau\) =\( 1+t\)^{-n/2}.
\end{equation*}
This case is similar to the second example, provided that we consider
positive times.

\section{On the hydrodynamic limit}
\label{sec:hydro}
In this paragraph, we consider the setting of either
Theorem~\ref{theo:main} or Theorem~\ref{theo:cubicquintic}. That is,
the semi-classical limit is justified for small time in Zhidkov
spaces.
Let $\Phi =\phi_{\rm eik}+\phi$, $\bv = \nabla \Phi$ and $\rho
=|a|^2$. As is easily checked, $(\rho,\bv)$ solves the following
compressible Euler equation:
\begin{equation}
  \label{eq:euler}
  \left\{
    \begin{aligned}
      &\d_t \rho + \DIV (\rho \bv)=0\quad ;\quad  
      \rho_{t=0}=|a_0|^2.\\
&\d_t \bv +\bv\cdot \nabla \bv + \nabla V +\nabla f\(\rho\) =0\quad ;\quad
      \bv_{\mid t=0}=\nabla \phi_0,
    \end{aligned}
\right.
\end{equation}
where $f(\rho)=\rho-1$ in the cubic case, and $f(\rho)=\rho^2+\l \rho$
in the cubic-quintic case.
To simplify the discussion, assume in this paragraph that $V_{\rm
  quad}=\phi_{\rm quad}=0$, hence $\phi_{\rm eik}=0$.
Proposition~\ref{prop:error1} implies in particular the convergence of
  the main two quadratic quantities, as $\eps \to 0$:
  \begin{itemize}
  \item Density: $|u^\eps|^2 \to \rho $ in $L^\infty([0,T_*]\times
  \R^n)$.
  \item Momentum: $\IM (\eps \overline u^\eps\nabla u^\eps)\to \rho
  \bv$ in $L^\infty([0,T_*]\times
  \R^n)$.
  \end{itemize}
It should be noted that if we assume only that for some $s>n/2+2$,
\begin{equation*}
  \|a_0^\eps -a_0\|_{X^s}=\delta^\eps_0=o(1)\quad \text{as }\eps \to 0,
\end{equation*}
the proof of Proposition~\ref{prop:error1} shows that we have:
\begin{equation*}
    \|\nabla(\phi^\eps-\phi)\|_{L^\infty([0,T_*];X^s)} +
    \|a^\eps-a\|_{L^\infty([0,T_*];X^s)} =\O \(\eps +\delta_0^\eps\).
  \end{equation*}
Therefore,
\begin{equation*}
  |u^\eps|^2 =\rho +\O \(\eps +\delta_0^\eps\)\quad ;\quad \IM (\eps
   \overline u^\eps\nabla
   u^\eps) = \rho \bv +\O \(\eps +\delta_0^\eps\).
\end{equation*}
To have a more precise asymptotics, it is necessary to work with the
assumption of Theorem~\ref{theo:main}. If for some $s>n/2+2$,
\begin{equation*}
  \|a_0^\eps -a_0-\eps a_1\|_{X^s}=\delta^\eps_1=o(\eps)\quad \text{as
  }\eps \to 0,
\end{equation*}
we get:
\begin{align*}
  |u^\eps|^2 &=\rho +2\eps\RE \(\overline a a^{(1)}\) +\O\(\eps^2 +
 \delta^\eps_1\) .\\
 \IM (\eps \overline u^\eps\nabla
   u^\eps) &= \rho \bv +\eps\(2  \RE \(\overline a a^{(1)}\) \bv +
  \rho \nabla \phi^{(1)}\) +\O\(\eps^2 +
 \delta^\eps_1\).
\end{align*}
Finally, note that in general, even if $a_1=0$, the modulation
$\phi^{(1)}$ is not trivial. Suppose that $a_1=0$: \eqref{eq:systlinear}
shows that $\d_t a^{(1)}_{\mid t=0}\not =0$, because of the source term
$\frac{i}{2}\Delta a$. Therefore,  even if
$\phi^{(1)}_{\mid t=0} = \d_t \phi^{(1)}_{\mid t=0}=0$, we have $\d_t^2
\phi^{(1)}_{\mid t=0}\not =0$ in general, and the correctors of
order $\eps$ in the above asymptotics are not trivial.
\smallbreak

However, if $a_0$ is \emph{real-valued} and $a_1=0$, then $a$ is
real-valued, $a^{(1)}$ is purely imaginary, so $\phi^{(1)}\equiv
0$. The same holds if $a_0$ is real-valued and $a_1$ is purely
imaginary.

\bigbreak

We end this section by studying the hydrodynamic limit in the case
when $\Omega\subset \R^n$ is a regular domain with bounded boundary
$\partial \Omega$ and $n\in \{2,3\}$ (either a bounded domain or an
exterior domain). 
To simplify the presentation, we consider the case without external potential
and without linear or quadratic initial phase.
The Gross--Pitaevskii equation is then supplemented
with the Neumann boundary condition:
\begin{equation}\label{GP2}
\left\{
\begin{aligned}
&i\eps \d_t u^\eps + \frac{\eps^2}{2}\Delta u^\eps =\( |u^\eps|^2-1\)u^\eps
&&\text{in }\Omega, \\
& \frac{\partial u^{\eps}}{\partial n}=0
&&\text{on }\partial \Omega, \\
\end{aligned}
\right.
\end{equation}
where $n$ is the unit outward normal to $\partial \Omega$.
Consider the corresponding limit system
\begin{equation}
\label{euler2}
\left\{
\begin{aligned}
&\d_t \rho + \DIV (\rho \nabla\phi)=0
&&\text{in }\Omega,\\
&\d_t \phi +\frac{1}{2}\left| \nabla \phi\right|^2 +\rho-1 =0
&&\text{in }\Omega, \\
&\nabla\phi\cdot n=0
&&\text{on }\partial \Omega.
\end{aligned}
\right.
\end{equation}

In~\cite{LinZhang}, Lin and Zhang proved that if $n=2$, then the
quadratic observables 
$|u^\eps|^2$ and $\eps\IM \(\overline u^\eps\nabla u^\eps\)$
converge towards the density $\rho$ and the momentum $\rho \nabla\phi$.
In the spirit of the pionnering work of Brenier \cite{Brenier},
in \cite{LinZhang} the strategy of the proof is to estimate
the modulated energy functional
$$
E^\eps \defn \frac{1}{\eps^2}\int_{\Omega}
\left| \eps\nabla u^\eps-i u^\eps \nabla \phi\right|^2
+ \(| u^\eps |^2-\rho\)^2\, dx.
$$
The assumption $n=2$ does not enter into the analysis of $E^\eps$
and only corresponds to the fact that they used
the Brezis--Gallou\"et inequality (see also \cite{ColinSoyeur})
to define sufficiently smooth solutions to the Gross--Pitaevskii equation.
They are now several 3D results (see \cite{PG05,Anton,Gallo,Gallo2}), and hence
one can justify the hydronamic limit for $n\in\{2,3\}$. In particular,
Theorem~\ref{theo:hydro} below is not new, but rather an update.
Yet, our main purpose here is to establish a local version of
the modulated energy functional. This is done in the proof
of Theorem~\ref{theo:hydro} (see~\eqref{eq:localenergy}),
by following the approach introduced in \cite{ACBKW}.

\begin{theorem}\label{theo:hydro}
Let $u^\eps$ and $(\rho,\phi)$
be classical solutions of \eqref{GP2} and \eqref{euler2} satisfying,
for some fixed $T>0$,
\begin{align*}
&u^\eps\in C([0,T];X^2(\Omega)),\quad |u^\eps|^2-1\in C([0,T];L^{2}(\Omega)),\\
&\rho\in C([0,T];X^{1}(\Omega)), \quad \rho-1\in C([0,T];L^{2}(\Omega)),\\
&\nabla\phi,\nabla^2\phi,\nabla^3\phi \in C([0,T];L^2(\Omega)\cap
L^\infty(\Omega)). 
\end{align*}
Assume that initially
$$
\left\| \eps\nabla u^\eps_{0}-i u^\eps_{0}\nabla\phi_{0}\right\|_{L^2(\Omega)}
+\left\| | u^\eps_{0}|^2-\rho_{0}\right\|_{L^2(\Omega)}=\mathcal{O}(\eps),
$$
then
\begin{align*}
&\left| u^{\eps}\right|^2-\rho =\mathcal{O}(\eps) \quad\text{in}\quad
L^{\infty}([0,T];L^{2}(\Omega)),\\ 
&\eps\IM (\overline u^\eps\nabla u^{\eps})-\rho
\nabla\phi=\mathcal{O}(\eps) \quad\text{in}\quad 
L^{\infty}([0,T];L^{1}_{\rm loc}(\Omega)).      
\end{align*}
\end{theorem}
\begin{proof}
The idea consists in filtering out the oscillations by the change of unknown
\begin{equation*}
 a^\eps(t,x)\defn u^\eps(t,x) e^{-i\phi(t,x)/\eps}.
\end{equation*}
The amplitude $ a^\eps$ solves
\begin{equation*}
 \d_t  a^\eps +\nabla\phi\cdot\nabla  a^\eps
+\frac{1}{2} a^\eps \Delta\phi- i\frac{\eps}{2}\Delta  a^\eps
=-\frac{i}{\eps}\( \left|  a^\eps \right|^{2} -
\rho\)  a^\eps .
\end{equation*}

Next set
$$
q^\eps \defn \frac{| a^\eps|^2-\rho}{\eps}\cdot
$$
We easily find that
\begin{equation*}
\d_{t} q^\eps + \DIV \(\IM \(\overline a^\eps\nabla a^\eps\)\) + \DIV
( q^\eps \nabla\phi)=0. 
\end{equation*}
Furthermore, with this notation, the equations for 
${\psi^\eps}\defn\nabla a^\eps$ read
\begin{equation*}
\d_t {\psi^\eps} +\nabla\phi\cdot\nabla {\psi^\eps}
+\frac{1}{2}{\psi^\eps} \Delta\phi +
{\psi^\eps}\cdot\nabla \nabla\phi+\frac{1}{2} a^\eps\nabla\Delta\phi  \\
+ i q^\eps {\psi^\eps} + i a^\eps \nabla
q^\eps=i\frac{\eps}{2}\Delta {\psi^\eps}.
\end{equation*}
Also, note that $\psi^{\eps}\cdot n=e^{-i\phi/\eps}(\nabla
u^\eps-i\eps^{-1}u^{\eps}\nabla\phi)\cdot n=0$ on $\partial\Omega$. 

We now introduce the modulated energy
$$
e^\eps  \defn | \psi^\eps|^2 + (q^\eps )^2.
$$
The key point is that
$$
q^\eps\DIV \big(\IM (\overline a^\eps\nabla a^\eps )\big)
+\RE \big( i a^\eps (\nabla q^\eps)\cdot\overline \psi^\eps\big)=
\DIV\big(\IM ( q^\eps \overline a^\eps {\psi^\eps} )\big).
$$
Hence, directly from the previous equations, we have
\begin{multline}\label{eq:localenergy}
\d_{t}e^\eps
+ \DIV(e^\eps\nabla\phi) + \DIV\bigl(2\IM ( q^\eps \overline a^\eps
{\psi^\eps} )\bigr) + \DIV\( \eps \IM \(\overline
\psi^\eps\cdot\nabla\psi^\eps\)\)   \\ 
=-( q^\eps)^2\Delta\phi-\RE\((
2{\psi^\eps}\cdot\nabla\nabla\phi+a^\eps\nabla\Delta\phi)\cdot\overline
\psi^\eps\). 
\end{multline}

We claim that
\begin{equation}\label{mei}
E^\eps(t)=\|e^\eps(t)\|_{L^1(\Omega)} \le
\|e^\eps(0)\|_{L^1(\Omega)}\exp\(Ct\)+C, 
\end{equation}
for some constant $C$ independent of $\eps$.
Since $v\cdot n=0$ and $\psi^\eps\cdot n=0$ on $\partial\Omega$,
by integrating in space and using the Gronwall's lemma,
to prove \eqref{claimhydro}, the only delicate point is to prove that,
\begin{equation}\label{claimhydro}
\int \left| a^\eps\nabla\Delta\phi\cdot\overline \psi^\eps \right|\,
dx \le C \| e^{\eps}\|_{L^1(\Omega)}+C. 
\end{equation}
To do so, as in Lemma~$1$ in \cite{PG05}, let $\chi\in C_{0}(\C)$ be
such that $0\le \chi\le 1$, 
$\chi(z)=1$ for $|z|\le 2$, and $\chi(z)=0$ for $|z|\ge 3$. Then write
$a^\eps=b^\eps+c^\eps$ where $b^\eps=\chi(a^\eps)a^\eps$ and
$c^{\eps}=(1-\chi(a^\eps))a^\eps$. We have 
$| b^\eps|\le 3$, $|c^\eps|\le \left| |a^\eps|^2-1\right|$ and hence
$$
\left\| b^\eps\right\|_{L^{\infty}(\Omega)}\le 3, \qquad
\left\| c^{\eps}\right\|_{L^{2}(\Omega)}\le
\left\| |a^\eps|^{2}-1\right\|_{L^{2}(\Omega)}
\le \eps\left\| q^\eps\right\|_{L^{2}(\Omega)}+\left\|
  \rho-1\right\|_{L^{2}(\Omega)}. 
$$
The desired estimate \eqref{claimhydro} then follows from
\begin{align*}
\left\| b^\eps\nabla\Delta\phi\cdot\overline \psi^\eps
\right\|_{L^1(\Omega)}\le  \left\| b^\eps\right\|_{L^{\infty}(\Omega)} 
\| \nabla\Delta\phi \|_{L^{2}(\Omega)}\left\|\psi^\eps
\right\|_{L^2(\Omega)},\\ 
\left\| c^\eps \nabla\Delta\phi \cdot\overline \psi^\eps \right\|_{L^1(\Omega)}
\le  \left\| c^\eps\right\|_{L^{2}(\Omega)} \| \nabla\Delta\phi
\|_{L^{\infty}(\Omega)} 
\left\|\psi^\eps \right\|_{L^2(\Omega)},
\end{align*}
and the elementary inequality $\sqrt{x}\le 1+x$.

Since
$$
e^\eps= \frac{1}{\eps^2}\left| \eps\nabla u^\eps-i
  u^\eps\nabla\phi\right|^2 + \frac{1}{\eps^{2}}\(| u^\eps
|^2-\rho\)^2, 
$$
the family $\(e^\eps(0)\)_{\eps\in ]0,1]}$ is bounded in $L^1(\Omega)$
by assumption. 
Consequently, it follows from \eqref{mei} that $\(e^\eps\)_{\eps\in ]0,1]}$
is bounded in $L^\infty( [0,T];L^1(\Omega))$.

By definition, this implies that $| u^\eps |^2-\rho=\mathcal{O}(\eps)$
in $L^\infty([0,T];L^{2}(\Omega))$. 
It remains to prove that 
$$
\eps\IM (\overline u^\eps\nabla u^{\eps})-\rho \nabla\phi 
=\mathcal{O}(\eps)\text{ in }L^{\infty}([0,T];L^{1}_{\rm loc}(\Omega)).
$$
Write
$$
\eps\IM (\overline u^\eps\nabla u^{\eps})-\rho \nabla\phi = \eps\IM
(\overline a^\eps\nabla a^{\eps}) 
+\(|a^\eps|^2-\rho\)\nabla\phi.
$$
Since $\nabla\phi\in L^\infty([0,T]\times\Omega)$, the previous result
implies that 
the second term is $\mathcal{O}(\eps)$ in $L^\infty([0,T];L^{2}(\Omega))$.
With regards to the first one, again write $a^\eps=b^\eps+c^\eps$
and use the obvious estimates 
\begin{align*}
&\| \eps\IM (\overline b^\eps\nabla a^{\eps}) \|_{L^{2}(\Omega)}\le
\eps \| b^\eps\|_{L^{\infty}(\Omega)} \|\nabla a^{\eps} \|_{L^{2}(\Omega)}
\le 3\eps \left\| e^\eps\right\|_{L^1(\Omega)}^{1/2},\\
&\| \eps\IM (\overline c^\eps\nabla a^{\eps}) \|_{L^{1}(\Omega)}\le
\eps \| c^\eps\|_{L^{2}(\Omega)} \|\nabla a^{\eps} \|_{L^{2}(\Omega)}
\le C\eps \left\| e^\eps\right\|_{L^1(\Omega)}^{1/2}+C\eps^2 \left\|
  e^\eps\right\|_{L^1(\Omega)}.
\end{align*}
This completes the proof.
\end{proof}

\section{Cubic-quintic nonlinearity}
\label{sec:quintic}

In view of Theorem~\ref{theo:cubicquintic}, we now consider
\eqref{eq:cubicquintic} in the case where the elliptic region becomes
relevant: $\l<0$, and assume for instance that there exists
$\underline x\in
\R^n$ such that $|a_0(\underline x)|^2<|\l|/2$. If we write $u^\eps
=a^\eps e^{i\Phi^\eps/\eps}$, where $(a^\eps,\Phi^\eps)$ is given by
\eqref{eq:systemmanuel}, then we naturally have to consider the limit
system:
\begin{equation}\label{eq:systlimcq}
  \left\{
\begin{aligned}
    \partial_t \phi +\frac{1}{2}\left|\nabla
    \phi\right|^2 +
    f_\l \(|a|^2\)= 0\quad &; \quad
    \phi\big|_{t=0}=\phi_0,\\
\partial_t a +\nabla \phi \cdot \nabla
    a +\frac{1}{2}a
\Delta \phi  = 0\quad & ;\quad
a\big|_{t=0}= a_0\, .
\end{aligned}
\right.
\end{equation}
Setting $v=\nabla \phi$, we find:
\begin{equation}\label{eq:systlimcq2}
  \left\{
\begin{aligned}
    \partial_t v +v\cdot \nabla v +
    \nabla f_\l \(|a|^2\)= 0\quad &; \quad
    v\big|_{t=0}=\nabla\phi_0,\\
\partial_t a +v \cdot \nabla
    a +\frac{1}{2}a
\DIV v  = 0\quad & ;\quad
a\big|_{t=0}= a_0\, .
\end{aligned}
\right.
\end{equation}
Then \cite[Theorem~3.2]{GuyCauchy} shows that \eqref{eq:systlimcq2} is
strongly ill-posed in Sobolev spaces. The problem remains in Zhidkov
spaces, since analyticity is essentially necessary. Indeed, Hadamard's
argument (see \cite{GuyCauchy} and references therein) shows for
instance that if $\phi_0$ is analytic near $\underline x$, then
\eqref{eq:systlimcq2}  has a $C^1$-solution only if $a_0$ is also
analytic near $\underline x$. So it may happen that
\eqref{eq:systlimcq2}  has no solution in $X^s$, even for $s$ large.
\smallbreak

On the other hand, if one is ready to work with analytic regularity,
then it becomes possible to justify the semi-classical limit for
\eqref{eq:cubicquintic}; see \cite{PGX93,ThomannAnalytic}.

\end{document}